\documentclass[11pt,reqno]{amsart}
\usepackage{fixltx2e}                    
\usepackage[osf]{mathpazo}  
\usepackage[utf8]{inputenc}            
\usepackage{amsmath}                     
\usepackage{amssymb, latexsym, stmaryrd, amsthm, dsfont, amsfonts, amsbsy,amsthm, amsmath, mathrsfs}            
\usepackage{mathtools}                   
\usepackage{bm}                          
\usepackage{enumerate}                   
\usepackage{verbatim}                    
\usepackage{url}   
\usepackage{lscape}                      

\usepackage{microtype}  
\usepackage[all, knot]{xy}
        \xyoption{arc} 
        \xyoption{web}                 
\makeatletter                            
\def\MT@register@subst@font{\MT@exp@one@n\MT@in@clist\font@name\MT@font@list
 \ifMT@inlist@\else\xdef\MT@font@list{\MT@font@list\font@name,}\fi}
\makeatother

\usepackage[pdftex,bookmarks,bookmarksnumbered,linktocpage,   
         colorlinks,linkcolor=blue,citecolor=blue]{hyperref}

\newcommand{\bit}{\begin{itemize}}    
\newcommand{\eit}{\end{itemize}}
\newcommand{\ben}{\begin{enumerate}}
\newcommand{\een}{\end{enumerate}}
\newcommand{\benormal}{\ben[\normalfont 1.]}   

\let\enormal\een
\newcommand{\benroman}{\ben[\normalfont (i)]}  
\let\eroman\een
\newcommand{\benbullet}{\ben[\textbullet]}     
\let\ebullet\een
\newcommand{\bde}{\begin{description}}
\newcommand{\ede}{\end{description}}

\newcommand{\?}{\ensuremath{\mkern0.4\thinmuskip}}   

\let\leq=\leqslant
\let\geq=\geqslant
\let\Box=\square                            

\let\epsilon=\varepsilon
\let\Lambda\varLambda
\let\Gamma\varGamma
\let\Delta\varDelta
\let\Lambda\varLambda
\let\Omega\varOmega
\let\Theta\varTheta
\let\Xi\varXi
\let\Pi\varPi
\let\Sigma\varSigma

\let\class=\mathsf                              
\let\oper=\mathbb                               

\bmdefine{\A}{A}                                
\bmdefine{\2}{2}
\bmdefine{\B}{B}
\bmdefine{\D}{D}
\bmdefine{\M}{M}                                
\bmdefine{\LLL}{L}                              
\bmdefine{\Fm}{Fm}                              
\bmdefine{\zerou}{[0{,}1]}  
\bmdefine{\T}{T}                                

\newcommand{\PSD}{\oper{P}_{\!\textsc{sd}}^{}}
\newcommand{\PPU}{\oper{P}_{\!\textsc{u}}^{}}
\newcommand{\PPR}{\oper{P}_{\!\textsc{r}}^{}}
\newcommand{\SSS}{\oper{S}}

\newcommand{\Mod}{\class{Mod}}

\bmdefine{\boldstar}{\mathchoice{\textstyle*}{\textstyle*}{\textstyle*}{\scriptstyle*}}
\newcommand{\Modstar}{\Mod^{\boldstar}}

\bmdefine{\btau}{\tau}                                  
\bmdefine{\brho}{\rho}                                  

\bmdefine{\leibniz}{\Omega}        

\bmdefine{\frege}{\Lambda}         

\makeatletter
\newcommand{\tarskidsp}{\mathord%
   {\m@th\raisebox{0pt}[0pt][0pt]{$\stackrel%
   {\raisebox{-2.7pt}[0ex][0pt]{$\displaystyle \,\?\thicksim$}}%
   {\displaystyle\leibniz}$}}}
\newcommand{\tarskitxt}{\mathord%
   {\m@th\raisebox{0pt}[0pt][0pt]{$\stackrel%
   {\raisebox{-2.7pt}[0ex][0pt]{$\,\?\thicksim$}}{\displaystyle\leibniz}$}}}
\newcommand{\tarskiscr}{\mathord%
   {{\m@th\raisebox{0pt}[0pt][0pt]{$\stackrel%
   {\raisebox{-2.4pt}[0ex][0pt]{$\scriptstyle \,\?\thicksim$}}%
   {\scriptstyle\leibniz}$}}}}
\newcommand{\tarskiscrscr}{\mathord%
   {{\m@th\raisebox{0pt}[0pt][0pt]{$\stackrel%
   {\raisebox{-2pt}[0ex][0pt]{$\scriptscriptstyle \,\?\thicksim$}}%
   {\scriptscriptstyle\leibniz}$}}}}
\newcommand{\tarski}{\@ifnextchar ^ %
   {\mathchoice{\tarskidsp\kern-.07em}{\tarskitxt\kern-.07em}%
   {\tarskiscr\kern-.07em}{\tarskiscrscr\kern-.07em}}%
   {\mathchoice{\tarskidsp}{\tarskitxt}{\tarskiscr}{\tarskiscrscr}}}
\makeatother

\theoremstyle{theorem}
\newtheorem{Theorem}{Theorem}[section]
\newtheorem{Lemma}[Theorem]{Lemma}
\newtheorem{Corollary}[Theorem]{Corollary}

\theoremstyle{definition}
\newtheorem{law}[Theorem]{Definition}
\newtheorem{exa}[Theorem]{Example}

\theoremstyle{remark}

\newtheorem{problem}{\bf Problem}

\newcommand{\C}{\boldsymbol{C}} 

\begin{document}
\title[On the complexity of the Leibniz hierarchy]{On the complexity of the Leibniz hierarchy}

\author{Tommaso Moraschini}
\address{Institute of Computer Science, Academy of Sciences of Czech Republic, Pod Vod\'arenskou v\v{e}\v{z}\'{i} $271/2$, $182$ $07$ Prague $8$, Czech Republic}
\email{moraschini@cs.cas.cz}
\date{\today}

\maketitle


\begin{abstract}
We prove that the problem of determining whether a finite logical matrix determines an algebraizable logic is complete for \textbf{EXPTIME}. The same result holds for the classes of order algebraizable, weakly algebraizable, equivalential and protoalgebraic logics. Finally, the same problem for the class of truth-equational logic is shown to be hard for \textbf{EXPTIME}.
\end{abstract}


\section{Introduction}

Abstract algebraic logic is a field that studies uniformly propositional logics \cite{Cz01,AAL-AIT-f,FJa09,FJaP03b}.\ One of its main achievements is the development of the so-called \textit{Leibniz hierarchy} (see Figure \ref{Hierarchy}), which provides a taxonomy of propositional logics inspired by the idea that logics can be classified accordingly to the definability of \textit{logical equivalence} and \textit{truth predicates}. The most prominent example of a class of logics belonging to the Leibniz hierarchy is probably the one of \textit{algebraizable logics} \cite{BP89}, i.e.\ logics which are equivalent to equational consequences relative to classes of algebras in the sense of \cite{BlJo99,BlJo06}.

The role of the Leibniz hierarchy in abstract algebraic logic has been compared \cite{JGRa11} to that of the \textit{Maltsev hierarchy} \cite{GaTa84,KeKi06,Tay73} in universal algebra, where the latter provides a taxonomy of varieties by means of properties typically related to the shape of congruence lattices. Recently, this analogy has been made precise in such a way that the Leibniz hierarchy can be regarded as a natural extension of the Maltsev hierarchy to arbitrary propositional logics \cite{JaMor19-1,JaMor19lh,JaMor19-3}. 

In the context of these considerations, the issue of whether it is actually possible to classify logics in the Leibniz hierarchy seems to deserve special attention.\ The goal of this paper is to shed light on this question, complementing previous work both in the setting of the Leibniz \cite{Mor16a} and Maltsev hierarchies  \cite{FV09} (see also \cite{Ho91c}).

It is well known that logics can be defined both syntactically, e.g.\ by means of Hilbert calculi, and semantically, e.g.\ by means of (logical) matrices.\ Accordingly, the  problem of classifying logics in the Leibniz hierarchy can be formulated in two versions: for each level $\class{K}$ of the Leibniz hierarchy of Figure \ref{Hierarchy} we consider the following problems:
\benormal
\item Is the logic determined by a given finite Hilbert calculus in $\mathsf{K}$?
\item Is the logic determined by a given finite matrix in $\mathsf{K}$?
\enormal
Problem 1 was shown to be undecidable in \cite{Mor16a}.\ On the other hand, Problem 2 turns out to solvable by means of algorithms that run in exponential or double exponential time (Lemmas \ref{Lem:EXPTIME1} and \ref{Lem:EXPTIME2}), see also \cite{TMo15}. Thus it makes sense to ask which is its \textit{computational complexity} \cite{ArBa09,Pa94}.

The main result of this paper shows that Problem 2 is hard for \textbf{EXPTIME} for every choice of $\mathsf{K}$ in Figure \ref{Hierarchy} (Theorem \ref{Thm:Hard}). Since the class \textbf{EXPTIME} is strictly larger than \textbf{PTIME} by the Hierarchy Theorem of Hartmanis and Stearns \cite{HaSt65}, this implies that these problems are not tractable.

\section{The Leibniz hierarchy}

This section contains a concise presentation of the main concepts used in the paper; for a systematic exposition we refer the reader to \cite{Cz01,AAL-AIT-f,FJa09,FJaP03b}. A \textit{matrix} is a pair $\langle \A, F\rangle$, where $\A$ is an algebra and $F \subseteq A$. A congruence $\theta$ of $\A$ is \textit{compatible} with $F \subseteq A$, when $F$ is a union of blocks of $\theta$. The largest congruence of $\A$ compatible with $F$ (which always exists) is denoted by $\leibniz^{\A}F$ and is called the \textit{Leibniz congruence} of $\A$ w.r.t.\ $F$. 
 A matrix $\langle \A, F\rangle$ is \textit{reduced} when $\leibniz^{\A}F$ is the identity relation. The \textit{reduction} of $\langle \A, F\rangle$ is the matrix $\langle \A / \leibniz^{\A}F, F / \leibniz^{\A}F\rangle$. The last matrix is always reduced. Given a class of matrices $\class{K}$, we denote by $\class{K}^{\boldstar}$ the class of isomorphic copies of the reductions of members of $\class{K}$, and by $\SSS(\class{K})$ the class of substructures of members of $\class{K}$.

A \textit{logic} $\vdash$ is a substitution invariant closure relation on the set of terms (constructed from a given countable set of variables) of a fixed algebraic language. A logic $\vdash$ is \textit{finitary} if for every set of formulas $\Gamma \cup \{ \varphi \}$,
\[
\Gamma \vdash \varphi \Longleftrightarrow \Delta \vdash \varphi\text{ for some finite }\Delta \subseteq \Gamma.
\]
Every class $\class{M}$ of matrices determines a logic $\vdash_{\class{M}}$ as follows:
\[
\Gamma \vdash_{\class{M}}\varphi \Longleftrightarrow \text{ for every evaluation }f,\text{ if }f[\Gamma] \subseteq F\text{, then }f(\varphi) \in  F.
\]
Observe that a matrix and its reduction determine the same logic. A matrix $\langle \A, F\rangle$ is a \textit{model} of a logic $\vdash$ when $\vdash \? \subseteq \? \vdash_{\langle \A, F\rangle}$. We denote by $\Modstar(\vdash)$ the class of reduced models of $\vdash$. 

The main classes in the \textit{Leibniz hierarchy} can de defined as follows. A logic $\vdash$ is \textit{protoalgebraic} if there is a set of formulas $\Delta(x, y, \vec{z})$ such that for every model $\langle \A, F\rangle$ of $\vdash$ and every $a, b \in A$ we have:
\begin{equation}\label{Eq:Protoalgebraic}
\langle a, b\rangle \in \leibniz^{\A}F \Longleftrightarrow \Delta(a, b, \vec{c}) \subseteq F \text{ for all }\vec{c} \in A.
\end{equation}
A logic is \textit{equivalential} if there exists a set of formulas $\Delta(x, y)$ without parameters $\vec{z}$, which satisfies condition (\ref{Eq:Protoalgebraic}).

\begin{figure}[t]
\[
\xymatrix@R=45pt @C=55pt @!0{
& {\txt{algebraizable}}\ar[dl]\ar[dr] & &  \\
{\txt{equivalential}}\ar[dr] & & {\txt{weakly\\ algebraizable}}\ar[dl]\ar[dr] & \\
 & {\txt{protoalgebraic}} & & {\txt{truth-equational}}
}
\]
\caption{The main classes in the Leibniz hierarchy.}
\label{Hierarchy}
\end{figure}
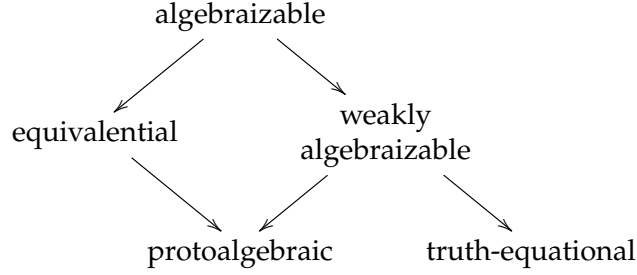

A logic $\vdash$ is \textit{truth equational} if there is a set of equations $\btau(x)$ such that for every $\langle \A, F\rangle \in \Modstar(\vdash)$ we have
\begin{equation}\label{Eq:TruthPredicates}
F = \{ a \in A : \A \vDash \btau(a) \}.
\end{equation}
Finally, a logic is \textit{algebraizable} (resp.\ \textit{weakly-algebraizable}) when it is equivalential (resp.\ protoalgebraic) and truth-equational. The Leibniz hierarchy, as described above, is represented in Figure \ref{Hierarchy} where arrows stand for the inclusion relation.

We will rely on the following observations.

\begin{Theorem}\label{Thm:Syntactic}
Let $\vdash$ be the logic determined by a finite set of finite matrices $\class{M}$.

\benormal
\item $\vdash$ is protoalgebraic if and only if there is a set of formulas $\Delta(x, y)$ (in variables $x$ and $y$) such that $\emptyset \vdash \Delta(x, x)$ and $x, \Delta(x, y) \vdash y$.
\item $\vdash$ is equivalential if and only if there is a set of formulas $\Delta(x, y)$ satisfying the conditions of point 1 and such for every basic $n$-ary operation $f$
\[
\Delta(x_{1}, y_{1}) \cup \dots \cup \Delta(x_{n}, y_{n}) \vdash \Delta(f(x_{1}, \dots, x_{n}), f(y_{1}, \dots, y_{n})).
\]
Moreover, in this case $\Delta$ satisfies condition \emph{(\ref{Eq:Protoalgebraic})} for every model of $\vdash$.
\item $\vdash$ is weakly algebraizable (resp.\ algebraizable) iff it is protoalgebraic (resp.\ equivalential) and there is a set of equations $\btau(x)$ such that for all $\langle \B, F \rangle \in \SSS(\class{M})$ and $b \in B$:
\[
b \in F \Longleftrightarrow \B / \leibniz^{\B}F \vDash \btau(b / \leibniz^{\B}F).
\] 
\enormal
\end{Theorem}

\begin{proof}
Items 1 and 2 are well-known (see for instance \cite{AAL-AIT-f}). Thus we detail only the proof of 3. First we claim that if $\vdash$ is protoalgebraic, then $\Modstar(\vdash) = \PSD ((\SSS(\class{M}))^{\boldstar})$. To prove this, observe that $\vdash$ is finitary logic, since $\class{M}$ is a finite set of finite matrices \cite[Theorem 4.4]{AAL-AIT-f}. In particular, this implies
\begin{equation}\label{Eq:model-theoretic-description-Modstar}
\Modstar(\vdash) = (\SSS\PPR(\class{M}))^{\boldstar} = (\PSD\SSS\PPU(\class{M}))^{\boldstar} = (\PSD\SSS(\class{M}))^{\boldstar},
\end{equation}
where $\PPR, \PSD$ and $\PPU$ are the class-operators of reduced, subdirect, and ultraproducts \cite[Theorem 4.7]{AAL-AIT-f}. Now, it is not hard to see that $(\PSD\SSS(\class{M}))^{\boldstar} \subseteq (\PSD((\SSS(\class{M}))^{\boldstar}))^{\boldstar} \subseteq \Modstar(\vdash)$. Hence we obtain 
\[
\Modstar(\vdash) = (\PSD((\SSS(\class{M}))^{\boldstar}))^{\boldstar}. 
\]
Recall that if $\vdash$ is protoalgebraic, then $\Modstar(\vdash)$ is closed under $\PSD$ \cite[Theorem 6.17]{AAL-AIT-f} (see also \cite{BP92,Cz80,DeJa96}). Together with the fact that $(\SSS(\class{M}))^{\boldstar} \subseteq \Modstar(\vdash)$, this implies that if $\vdash$ is protoalgebraic, then
\[
\Modstar(\vdash) = (\PSD((\SSS(\class{M}))^{\boldstar}))^{\boldstar} = \PSD((\SSS(\class{M}))^{\boldstar}),
\]
establishing the claim.

From the claim it follows that if $\vdash$ is protoalgebraic, then $\vdash$ is truth-equational if and only if there is a set of equations $\btau(x)$ such that
for all $\langle \B, F \rangle \in \SSS(\class{M})$ and $b \in B$:
\[
b \in F \Longleftrightarrow \B / \leibniz^{\B}F \vDash \btau(b / \leibniz^{\B}F).
\] 
This easily implies that condition 3 holds.
\end{proof}

\begin{Lemma}\label{Lem:TestTruthEquationality}
Let $\vdash$ be the logic determined by a finite set of finite matrices $\class{M}$. Then $\vdash$ is truth-equational if and only if there is a set of equations $\btau(x)$ such that for every $\langle \B, G \rangle \in \Modstar(\vdash)$, where $\B$ is one-generated, and every $\langle \A, F \rangle \in \class{M}$ we have
\[
G = \{ b \in B : \B \vDash \btau(b) \} \text{ and }F \subseteq \{ a \in A : \A / \leibniz^{\A}F \vDash \btau(a / \leibniz^{\A}F) \}.
\]
In this case, we can take $\btau(x)$ to be the set of equations $\epsilon \thickapprox \delta$ such that
\[
\epsilon^{\A / \leibniz^{\A}F}(a/ \leibniz^{\A}F) = \delta^{\A/ \leibniz^{\A}F}(a/ \leibniz^{\A}F), \text{ for all }a\in F\text{ and } \langle \A, F \rangle \in \class{M}.
\]
\end{Lemma}

\begin{proof}
We claim that if $\epsilon(x) \thickapprox \delta(x)$ is an equation such that $\epsilon^{\A / \leibniz^{\A}F}(a / \leibniz^{\A}F) = \delta^{\A / \leibniz^{\A}F}(a / \leibniz^{\A}F)$ for all $a\in F$ and $\langle \A, F \rangle \in \class{M}$, then
\[
\epsilon^{\B}(b) = \delta^{\B}(b) \text{ for every }\langle \B, G \rangle \in \Modstar(\vdash) \text{ and }b \in G.
\]

To prove this, observe that $\Modstar(\vdash)$ can be seen as a class of first-order structures in the algebraic language of $\vdash$ extended with a unary predicate symbol $P(x)$. More precisely, a matrix $\langle \B, G \rangle \in \Modstar(\vdash)$ can be regarded as a first-order structure given by the algebra $\B$, equipped with the interpretation of $P(x)$ given by the set $G \subseteq B$. 

Now, consider a matrix $\langle \B, G \rangle \in \Modstar(\vdash)$, and an equation $\epsilon \thickapprox \delta$ such that $\epsilon^{\A / \leibniz^{\A}F}(a / \leibniz^{\A}F) = \delta^{\A / \leibniz^{\A}F}(a / \leibniz^{\A}F)$ for all $a\in F$ and $\langle \A, F \rangle \in \class{M}$. We have
\begin{equation}\label{Eq:sentence}
\class{M} \vDash \forall x (P(x) \to \epsilon(x) \thickapprox \delta(x) ).
\end{equation}
Moreover, as in (\ref{Eq:model-theoretic-description-Modstar}), we have
\begin{align*}
\Modstar(\vdash) &= (\PSD\SSS(\class{M}))^{\boldstar}.
\end{align*}
Hence $\langle \B, G \rangle$ can be identified with a matrix of the form $\langle \C / \leibniz^{\C}H, H / \leibniz^{\C}H \rangle$ for some $\langle \C, H \rangle \in \PSD\SSS (\class{M})$. 

Then consider an arbitrary element $c \in C$ such that $c / \leibniz^{\C}H \in H / \leibniz^{\C}H$. Since the sentence in (\ref{Eq:sentence}) is preserved by $\PSD$ and $\SSS$, it holds in $\langle \C, H \rangle$ as well. Together with the fact that $c \in H$, this implies that $\epsilon^{\C}(c) = \delta^{\C}(c)$. Hence we conclude that $\C/ \leibniz^{\C}H \vDash \epsilon(c/ \leibniz^{\C}H ) \thickapprox \delta(c/ \leibniz^{\C}H )$, establishing the claim.


Now, we move to the proof of the main statement. Indeed, it is enough to prove the ``if'' part. To this end, suppose that there is a set of equations $\btau(x)$ satisfying the assumptions. Then consider an arbitrary matrix $\langle \B, G \rangle \in \Modstar(\vdash)$. In order to prove that $\vdash$ is truth-equational, it will be enough to show that
\[
G = \{ b \in B : \B \vDash \btau(b) \}.
\] 
From the assumption we know that $F \subseteq \{ a \in A : \A / \leibniz^{\A}F \vDash \btau(a / \leibniz^{\A}F) \}$ for every $\langle \A, F \rangle \in \class{M}$. Together with the claim, this implies that $G \subseteq \{ b \in B \colon \B \vDash \btau(b) \}$. Conversely, consider $b \in B$ such that $\B \vDash \btau(b)$. Then consider the subalgebra $\C$ of $\B$ generated by $b$. As the congruence $(C \times C) \cap \leibniz^{\B}G$ is compatible with $G \cap C$, we obtain
\[
(C \times C) \cap \leibniz^{\B}G \subseteq \leibniz^{\C}(G \cap C).
\]
Together with the fact that $\B \vDash \btau(b)$, this implies
\begin{equation}\label{Eq:QuotientOfC}
C / \leibniz^{\C}(G \cap C) \vDash \btau(b / \leibniz^{\C}(G \cap C)).
\end{equation}
Now, observe that $\langle \C / \leibniz^{\C}(G \cap C), (G \cap C) / \leibniz^{\C}(G \cap C) \rangle$ is a reduced model of $\vdash$, and that $\C / \leibniz^{\C}(G \cap C)$ is one-generated. Hence we can apply the assumption to (\ref{Eq:QuotientOfC}), obtaining
\[
b/ \leibniz^{\C}(G \cap C) \in (G \cap C) / \leibniz^{\C}(G \cap C).
\]
Since $\leibniz^{\C}(G \cap C)$ is compatible with $G \cap C$, we get that $b \in G \cap C \subseteq G$. Hence we conclude that 
\[
G \supseteq \{ b \in B : \B \vDash \btau(b) \}
\]
and, therefore, that $\vdash$ is truth equational.

To conclude the proof, suppose that there is a set of equations $\btau(x)$ satisfying the condition in the statement. Moreover, let $\btau'(x)$ be the set of equations $\epsilon \thickapprox \delta$ such that $\epsilon^{\A / \leibniz^{\A}F}(a/ \leibniz^{\A}F) = \delta^{\A/ \leibniz^{\A}F}(a/ \leibniz^{\A}F)$, for all $a\in F$ and $\langle \A, F \rangle \in \class{M}$. Since $F \subseteq \{ a \in A \colon \A / \leibniz^{\A}F \vDash \btau(a\leibniz^{\A}F) \}$ for all $\langle \A, F \rangle \in \class{M}$, we obtain that $\btau(x) \subseteq \btau'(x)$. In particular, this implies that for every $\langle \B, G \rangle \in \Modstar(\vdash)$ such that $\B$ is one-generated,
\[
G \supseteq \{ b \in B : \B \vDash \btau'(b) \}.
\]
The reverse inclusion follows from the claim at the beginning of the proof.
\end{proof}

We will also make use of the following technical lemma:\footnote{I am indebted to the referee of my \cite{Mor16a} for suggesting a simplification in the proof of Lemma \ref{Lem:TestTruthEquationality}.}

\begin{Lemma}\label{Lem:TrivialCases}
Let $\vdash$ be the logic determined by a reduced matrix $\langle \A, F \rangle$.
\benormal
\item If $\A$ is trivial, then $\vdash$ is weakly algebraizable if and only if $F = A$.
\item If $\A$ is non-trivial and all its basic operations are constants, then $\vdash$ is not weakly algebraizable.
\item If $\A$ is trivial, then $\vdash$ is truth-equational if and only if $F = A$.
\item If $\A$ is non-trivial and all its basic operations are constants, then $\vdash$ is truth-equational if and only if there is a constant ${ \bf a}$ which is interpreted inside $F$.
\enormal
\end{Lemma}

\begin{proof}
1. Observe that if $\A$ is trivial, then $\vdash$ is either the \textit{inconsistent} logic, i.e. the logic whose set of tautologies coincides with the set of all formulas, or the \textit{almost inconsistent} logic, i.e. the logic with empty set of tautologies but such that $\varphi \vdash \psi$ for all formulas $\varphi$ and $\psi$. More in detail, $\vdash$ is inconsistent when $F  = A$, and it is almost inconsistent when $F = \emptyset$. From condition 3 of Theorem \ref{Thm:Syntactic} it follows that the inconsistent logic is weakly algebraizable, while the almost inconsistent one is not.

2.  Suppose that $\A$ is non-trivial and that all its basic operations are constants. Since all operations of $\A$ are constants, the Leibniz congruence of $\langle \A, F\rangle$ has only two blocks, i.e. $F$ and $A \smallsetminus F$. Together with the fact that $\langle\A, F\rangle$ is reduced and that $\A$ non-trivial, this implies that $\A$ is a two element algebra, say with universe $\{ 0, 1\}$, and $F = \{ 1 \}$. Hence $\vdash$ is determined by $\langle \A, \{ 1 \} \rangle$.

Suppose towards a contradiction that $\vdash$ is weakly algebraizable.\ Then there is a set of formulas $\Delta(x, y)$ as in condition 1 of Theorem \ref{Thm:Syntactic}.\ In particular, this means that
\[
\emptyset \vdash \Delta(x, x) \text{ and }x, \Delta(x, y) \vdash y.
\]
Due to the poorness of the language of $\A$, we have that $\Delta(x, y)$ can contain only constants ${ \bf a}, {\bf b}, {\bf c}\dots$ and the variables $x, y$. Now, since $\Delta(x, x)$ is a set of tautologies of $\vdash$ and $\vdash$ is determined by $\langle \A, \{ 1 \} \rangle$, we have that ${\bf a } \notin\Delta(x, y)$ for each  constant ${\bf a}$ which is interpreted to $0$. Similarly, as $x$ and $y$ can be interpreted as $0$ and $\Delta(x, x)$ is a set of tautologies, we obtain $x, y \notin \Delta(x, y)$. Then $\Delta(x, y)$ is a set of constants that are interpreted to $1$. Hence $\Delta(x, y)$ is a set of tautologies of $\vdash$, since $\vdash$ is determined by $\langle \A, \{ 1 \} \rangle$. This means that the rule $x, \Delta(x, y) \vdash y$ specializes to $x \vdash y$, which is easily seen to fail in $\langle \A, \{ 1 \} \rangle$.

3. Suppose that $\A$ is trivial. If $F =A$, then all reduced models of $\vdash$ are isomorphic to $\langle \A, F\rangle$. Hence $\vdash$ is truth-equational with $\btau=\emptyset$. Then suppose that $F = \emptyset$. Observe that there cannot be a set of equations $\btau$ whose set of solutions in $\A$ is empty. Hence $\vdash$ is not truth-equational.

4. Suppose that $\A$ is non-trivial. As motivated in the proof of condition 2, $\A$ is a two element algebra, say with universe $\{ 0, 1 \}$, and $F = \{ 1 \}$. If no constant symbol of $\A$ is interpreted to $1$, then $F$ cannot be defined as in (\ref{Eq:TruthPredicates}). Hence $\vdash$ is not truth-equational. On the other hand, if there is a constant ${ \bf a}$ which is interpreted to $1$, then $\Modstar(\vdash)$ is the class of isomorphic copies of $\langle \A, F\rangle$ and of the matrix $\langle \B, B\rangle$, where $\B$ is the trivial algebra. This means that $\vdash$ is truth-equational with $\btau(x) \coloneqq \{ x \thickapprox { \bf a} \}$.
\end{proof}

Recall that if $\A$ is an algebra, the free $n$-generated algebra $\boldsymbol{Tm}(x_{1}, \dots, x_{n})$ in the variety generated by $\A$ is isomorphic to the subalgebra of $\A^{(A^{n})}$ generated by the projections (see \cite{Be05k} if necessary). Sometimes we tacitly identify the elements of the free algebra $\boldsymbol{Tm}(x_{1}, \dots, x_{n})$ (i.e.\ the equivalence classes of formulas with variables among $x_{1}, \dots, x_{n}$ that have the same interpretation in $\A$) with particular formulas $\varphi(x_{1}, \dots, x_{n})$ that represent them.

\begin{Lemma}\label{Lem:concrete-sets-of-formulas}
Let $\langle \A, F \rangle$ be a finite reduced matrix, and $\boldsymbol{Tm}(x)$, $\boldsymbol{Tm}(x, y)$ respectively the free one-generated and the free two-generated algebras over the variety generated by $\A$. Consider also the following sets of formulas and equations
\begin{align*}
\Delta(x, y)  \coloneqq & \{ \varphi \in Tm(x, y) \colon \varphi^{\A}(a, a) \in F, \text{ for all }a \in A\}\\
\btau(x) \coloneqq & \{ \epsilon \thickapprox \delta \colon \epsilon, \delta \in Tm(x) \text{ and }\epsilon^{\A}(a) = \delta^{\A}(a), \text{ for all }a \in F\}.
\end{align*}
The following conditions hold for the logic $\vdash$ determined by $\langle \A, F \rangle$.
\benormal
\item $\vdash$ is protoalgebraic iff $\Delta$ satisfies the requirement in condition 1 of Theorem \ref{Thm:Syntactic}.
\item $\vdash$ is equivalential iff  $\Delta$ satisfies the requirement in condition 2 of Theorem \ref{Thm:Syntactic}.
\item $\vdash$ is weakly algebraizable (resp.\ algebraizable) iff  $\vdash$ is protoalgebraic (resp.\ equivalential) and $\btau$ satisfies the requirement in condition 3 of Theorem \ref{Thm:Syntactic}.
\item $\vdash$ is truth-equational iff $\btau$ satisfies the requirement in Lemma \ref{Lem:TestTruthEquationality}.
\enormal
\end{Lemma}

\begin{proof}
1. We have to prove that $\vdash$ is protoalgebraic if and only if $\emptyset \vdash \Delta(x, x)$ and $x, \Delta(x, y) \vdash y$. The ``if'' part is a consequence of Theorem \ref{Thm:Syntactic}(1). To prove the ``only if'' part, suppose that $\vdash$ is protoalgebraic. By condition 1 of Theorem \ref{Thm:Syntactic} there is a set of formulas $\Delta'(x, y)$ such that $\emptyset \vdash \Delta'(x, x)$ and $x, \Delta'(x, y) \vdash y$. Since the logic $\vdash$ is determined by $\langle \A, F \rangle$, we can assume w.l.o.g.\ that $\Delta' \subseteq Tm(x, y)$. Together with $\emptyset \vdash \Delta'(x, x)$ and the definition of $\Delta$, this implies $\Delta' \subseteq \Delta$. Bearing in mind that $x, \Delta'(x, y) \vdash y$, we conclude that $x, \Delta(x, y) \vdash y$. Finally, the fact that $\emptyset \vdash \Delta(x, x)$ is an immediate consequence of the definition of $\Delta$.

2. The ``if'' part of the statement is a consequence of Theorem \ref{Thm:Syntactic}(2). To prove the ``only if'' part, suppose that $\vdash$ is equivalential. By condition 2 of Theorem \ref{Thm:Syntactic} there is a set $\Delta'(x, y)$ such that
\begin{align*}
\emptyset &\vdash \Delta'(x, x) \qquad x, \Delta'(x, y) \vdash y\\
\Delta'(x_{1}, y_{1}) \cup \dots \cup \Delta&(x_{n}, y_{n}) \vdash \Delta'(f(x_{1}, \dots, x_{n}), f(y_{1}, \dots, y_{n}))
\end{align*}
for every basic $n$-ary operation $f$. To conclude the proof, it will be enough to show that the sets $\Delta(x, y)$ and $\Delta'(x, y)$ are interderivable in $\vdash$. 

As in the proof of 1, we can assume w.l.o.g.\ that $\Delta' \subseteq \Delta$. Therefore it only remains to show that $\Delta'(x, y) \vdash \varphi$ for every $\varphi \in \Delta$. To this end, consider $a, b \in A$ such that $\Delta'^{\A}(a, b) \subseteq F$. By the last part of condition 2 of Theorem \ref{Thm:Syntactic}, this implies $\langle a, b \rangle \in \leibniz^{\A}F$. Since $\langle \A, F \rangle$ is reduced, we conclude that $a = b$. Together with the definition of $\Delta$, this implies $\Delta^{\A}(a, b) = \Delta^{\A}(a, a) \subseteq F$. Since $\vdash$ is the logic determined by $\langle \A, F \rangle$, we conclude that $\Delta'(x, y) \vdash \varphi$ for all $\varphi \in \Delta(x, y)$.

4. The ``if'' part of the statement is a consequence of Lemma \ref{Lem:TestTruthEquationality}. To prove the ``only if'' part, suppose that $\vdash$ is truth-equational. From the definition of $\btau$ and the fact that $\langle \A, F \rangle$ is reduced it follows that
\[
F \subseteq	\{ a \in A \colon \A \vDash \btau(a) \} =  \{ a \in A \colon \A / \leibniz^{\A}F \vDash \btau(a/\leibniz^{\A}F) \}.
\]
Then consider $\langle \B, G \rangle \in \Modstar(\vdash)$ such that $\B$ is one-generated. It only remains to prove that for every $b \in B$,
\begin{equation}\label{Eq:selecting-equ1}
b \in G \Longleftrightarrow \B \vDash \btau(b).
\end{equation}

To this end, consider a matrix $\langle \B, G \rangle \in \Modstar(\vdash)$ such that $\B$ is one-generated. Moreover, let $\btau'$ be the set of equations $\epsilon(x) \thickapprox \delta(x)$ such that $\epsilon^{\A}(a) = \delta^{\A}(a)$ for every $a \in F$. Since $\vdash$ is truth-equational, we can apply the last part of the statement of Lemma \ref{Lem:TestTruthEquationality} and the fact that $\langle \A, F\rangle$ is reduced obtaining that for every $b \in B$,
\begin{equation}\label{Eq:selecting-equ2}
b \in G \Longleftrightarrow \B \vDash \btau'(b).
\end{equation}
Now, as in (\ref{Eq:model-theoretic-description-Modstar}), $\B$ belongs to the variety generated by $\A$. As a consequence,
\[
\B \vDash \btau(b) \Longleftrightarrow \B \vDash \btau'(b)
\]
for every $b \in B$. Together with (\ref{Eq:selecting-equ2}), this implies (\ref{Eq:selecting-equ1}) as desired.

3. We detail the case of weakly algebraizable logics only. The ``if'' part of the statement is a consequence of Theorem \ref{Thm:Syntactic}(3). To prove the ``only if'' part, suppose that $\vdash$ is weakly algebraizable. Clearly, $\vdash$ is protoalgebraic. Moreover, since $\vdash$ is truth-equational, we can apply point 4 yielding that $\btau$ satisfies the conditions in Lemma \ref{Lem:TestTruthEquationality}. This clearly implies that $\btau$ satisfies the (weaker) requirement in Theorem \ref{Thm:Syntactic}(3).
\end{proof}

\section{Upper bounds}

We begin by establishing some upper bound to the computational complexity of the problem of classifying logics determined by a finite reduced matrix of finite type into the Leibniz hierarchy of Figure \ref{Hierarchy}. To this end, let us explain how we represent the inputs of this problem, that is how we represent finite matrices $\langle \A, F \rangle$. The finite set $A$ could be represented by its cardinality, say $n$ (and $A$ can be taken to be $\{ 1, \dots, n \}$). A $k$-ary operation can be represented as a $k$-dimensional array of elements of $A$, i.e.\ a sequence of $n^{k}$ elements of $A$. Finally, $F$ can be represent by a natural number $m \leq n$, indicating that $F = \{ k \in \omega : 1 \leq k \leq m \}$ (observe that $F = \emptyset$, when $m = 0$).


The next result identifies upper bounds to the complexity of the problem of determining whether the logic of a finite reduced matrix is algebraizable:

\begin{Lemma}\label{Lem:EXPTIME1}
The problem of determining whether the logic of a finite reduced matrix of finite type is algebraizable (resp.\ equivalential, protoalgebraic, weakly algebraizable) is in \textbf{EXPTIME}.
\end{Lemma}

\begin{proof}
We detail the proof for the case of weakly algebraizable logics. Consider a finite reduced matrix $\langle \A, F\rangle$ of finite type, and let $\vdash$ be the logic it defines. Moreover, let $\Delta$ and $\btau$ be the sets defined in Lemma \ref{Lem:concrete-sets-of-formulas}. In order to determine whether $\vdash$ is weakly algebraizable or not, it is enough to check whether $\Delta$ and $\btau$ satisfy, respectively, the requirement in condition 1 and 3 of Theorem \ref{Thm:Syntactic} (see conditions 1 and 3 of Lemma \ref{Lem:concrete-sets-of-formulas} if necessary). This task can be done mechanically, since the free algebras $\boldsymbol{Tm}(x)$ and $\boldsymbol{Tm}(x, y)$ in the variety generated by $\A$ are finite and of size at most $\vert A^{A^{2}} \vert$.

It only remains to show that this test can be carried on in exponential time. To this end, let us describe it in more detail. The algorithm starts reading the input $\langle \A, F\rangle$. If it reads that $\A$ is trivial, then it checks whether $F$ is empty or not and provides an output according to point 1 of Lemma \ref{Lem:TrivialCases}. If $\A$ is non-trivial, then the algorithm starts checking whether $\A$ has a basic operation, which is not a constant. If this is not the case, then the algorithm establishes that $\vdash$ is not weakly algebraizable according to point 2 of Lemma \ref{Lem:TrivialCases}. This process takes time $O(n)$, where $n$ is the length of the input.

If the algorithm does not have halted yet, then $\A$ is non-trivial and has a basic operation, which is not a constant. In particular, this means that the length $n$ of the input dominates the cardinality of $\A$. First we need to construct the sets $Tm(x)$ and $Tm(x, y)$ respectively of all unary and binary term-functions of $\A$. This task amounts to that of constructing the one-generated and the two-generated free algebras of the variety generated by $\A$ (which are isomorphic to subalgebras of $\A^{A}$ and $\A^{(A^{2})}$, respectively). Since $n$ dominates the cardinality of $\A$, the time needed in the construction of $Tm(x)$ and $Tm(x, y)$ can be bounded exponentially in $n$. Thus the sets $\Delta$ and $\btau$ can be constructed in exponential time in $n$.

By Lemma \ref{Lem:concrete-sets-of-formulas}(1), $\vdash$ is protoalgebraic if and only if $\Delta$ satisfies the right-hand side of condition 1 of Theorem \ref{Thm:Syntactic}. Then the algorithm checks if this is the case or not (which can be done in exponential time). Accordingly, if $\vdash$ is not protoalgebraic, the algorithm stops and provides a negative answer.

If $\vdash$ is protoalgebraic, then by Lemma \ref{Lem:concrete-sets-of-formulas}(3) the logic $\vdash$ is weakly algebraizable if and only if $\btau$ satisfies the requirement in the display of condition 3 of Theorem \ref{Thm:Syntactic}. To this end, we construct all substructures $\langle \B, G \rangle$ of $\langle \A, F \rangle$ (which can be done in exponential time). For each such submatrix $\langle \B, G \rangle$, the Leibniz congruence $\leibniz^{\B}G$ can be computed in $n \log n$ time.\footnote{This is justified as follows. First observe that the matrix $\langle \B, G \rangle$ can be transformed into an automaton $\langle \B^{+}, G \rangle$, where $\B^{+}$ is the algebra with universe $B$ and whose basic operations are the unary functions of the form $f^{\B}(b_{1}, \dots, b_{n}, x, b_{n+1}, \dots, b_{m})$ for some $\vec{b} \in B$ and some basic operation $f$ of $\B$. Observe that the length of $\langle \B^{+}, G \rangle$ is essentially the same of $\langle \B, G \rangle$, i.e.\ $n$. Now, it is known that the task of computing the Leibniz congruence $\leibniz^{\B}G$ amounts to that of minimizing the number of states of the automaton $\langle \B^{+}, G \rangle$. Since the time needed to solve the latter task is bounded above by $n \log n$ \cite{Fr08f}, we are done.} Finally, for every $\langle \B, G \rangle$, checking the requirement in the display of condition 3 of Theorem \ref{Thm:Syntactic} is an easy task. If this requirement holds, then the algorithm provides a positive answer, otherwise it provides a negative one.

A similar argument shows that the problem of determining whether the logic of a finite reduced matrix of finite type is equivalential (resp.\ protoalgebraic, algebraizable) belongs to \textbf{EXPTIME}. 
\end{proof}

The next result identifies upper bounds to the complexity of determining whether the logic of a finite reduced matrix is truth-equational:

\begin{Lemma}\label{Lem:EXPTIME2}
The problem of determining whether the logic of a finite reduced matrix of finite type is truth-equational is in ${\bf 2}$\textbf{-EXPTIME}.
\end{Lemma}

\begin{proof}
Let  $\langle \A, F\rangle$ be a finite reduced matrix of finite type, and let $\vdash$ be the logic it determines. Let also $\btau$ be the set of equations in the statement of Lemma \ref{Lem:concrete-sets-of-formulas}.  In order to determine whether $\vdash$ is truth-equational or not, it is enough to check whether $\btau$ satisfies the requirement in Lemma \ref{Lem:TestTruthEquationality} (see conditions 4 of Lemma \ref{Lem:concrete-sets-of-formulas} if necessary). 

We proceed to explain why this task can be done mechanically. First observe that the free one-generated algebra $\boldsymbol{Tm}(x)$ in the variety generated by $\A$, and the set $\btau$ can be constructed mechanically as in the proof of Lemma \ref{Lem:EXPTIME1}. Then we construct all sets $\Gamma \subseteq Tm(x)$ such that $\langle \boldsymbol{Tm}(x), \Gamma \rangle$ is a model of $\vdash$. Then, we check whether for every such $\Gamma$ we have that
\begin{equation}\label{Eq:Algorithm-Truth-Eq}
\Gamma / \leibniz^{\boldsymbol{Tm}(x)}\Gamma = \{\varphi / \leibniz^{\boldsymbol{Tm}(x)}\Gamma \in Tm(x) / \leibniz^{\boldsymbol{Tm}(x)}\Gamma \colon  \btau(\varphi) \subseteq \leibniz^{\boldsymbol{Tm}(x)}\Gamma \}.
\end{equation}
Observe that the above equality can be checked mechanically, since $\boldsymbol{Tm}(x)$ is a finite algebra. The test says that $\vdash$ is truth-equational, if the above condition holds, and it says that $\vdash$ is not truth-equational otherwise.

The fact that this test works as expected is guaranteed by the following argument. As remarked in the proof of Lemma \ref{Lem:TestTruthEquationality}, we have that
\begin{align*}
\Modstar(\vdash) &= (\PSD\SSS(\langle \A, F \rangle))^{\boldstar}.
\end{align*} 
In particular, this means that the algebraic reducts of $\Modstar(\vdash)$ belong to the variety generated by $\A$. This observations has two consequences. On the one hand, it means that the free one-generated algebra of the variety generated by the algebraic reducts of $\Modstar(\vdash)$ is indeed $\boldsymbol{Tm}(x)$. On the other hand, it implies that the reduced models $\langle \B, G \rangle$ of $\vdash$ such that $\B$ is one-generated are (up to isomorphism) the ones of the form $\langle \boldsymbol{Tm}(x)/ \leibniz^{\boldsymbol{Tm}(x)}\Gamma, \Gamma / \leibniz^{\boldsymbol{Tm}(x)}\Gamma \rangle$ where $\Gamma \subseteq Tm(x)$ is such that $\langle \boldsymbol{Tm}(x), \Gamma \rangle$ is a model of $\vdash$. Hence our test amounts to checking whether $\btau$ satisfies the requirement in the statement of Lemma \ref{Lem:TestTruthEquationality}.

%
%


It only remains to show that this test can be carried on in double exponential time. To this end, let us describe it in more detail. The algorithm starts reading the input $\langle \A, F\rangle$. If it reads that $\A$ is trivial, then it checks whether $F$ is empty or not, and provides an output according to condition 3 of Lemma \ref{Lem:TrivialCases}. If $\A$ is non-trivial, then the algorithm starts checking whether all a basic operations of $\A$ are constants. If this is the case, then the algorithm checks whether at least one of them is interpreted inside $F$ and provides an output according to condition 4 of Lemma \ref{Lem:TrivialCases}. This process takes time $O(n)$, where $n$ is the length of the input $\langle \A, F\rangle$.

If the algorithm does not have halted yet, then $\A$ is non-trivial and has a basic operation, which is not a constant. In particular, this means that the length $n$ of the input dominates the cardinality of $\A$. From now on we just apply the algorithm described at the beginning of the proof. First observe that $\btau$ and $\boldsymbol{Tm}(x)$ can be constructed in exponential time in $n$, as in the proof of Lemma \ref{Lem:EXPTIME1}.

Then we need to construct all sets $\Gamma \subseteq Tm(x)$ such that   $\langle \boldsymbol{Tm}(x), \Gamma \rangle$ is a model of $\vdash$. This computation requires a number of constructions of order
\[
2^{n^{n}} \cdot n^{n} \cdot n \leq 2^{2^{n^{2} +2n+n}}
\]
since the number of  subsets $\Gamma \subseteq Tm(x)$ is bounded above by $2^{n^{n}}$, and for each of $\Gamma \subseteq Tm(x)$ we have to check that every formula $\varphi \in Tm(x)$ (there are at most $n^{n}$ of them) such that $\Gamma \vdash \varphi$ belongs to $\Gamma$ (checking whether $\Gamma \vdash \varphi$ involves at most $n$ assignments of the variable $x$ in $\A$, since $n$ dominates the cardinality of $A$). Since all these constructions can be carried on in polynomial time in the length $n$ of the input, we conclude that constructing all the sets $\Gamma \subseteq Tm(x)$ such that   $\langle \boldsymbol{Tm}(x), \Gamma \rangle$ is a model of $\vdash$ takes at most double exponential time in $n$.

Let $\langle \boldsymbol{Tm}(x), \Gamma \rangle$ be a model of $\vdash$. As shown in the proof of Lemma \ref{Lem:EXPTIME1}, the Leibniz congruence $\leibniz^{\boldsymbol{Tm}(x)}\Gamma$ can be computed in $m \log m$ time in the length of $\leibniz^{\boldsymbol{Tm}(x)}\Gamma$, which is bounded by $2^{p(n)}$ for some polynomial $p(n)$. Bearing this in mind, the time needed in the construction of $\leibniz^{\boldsymbol{Tm}(x)}\Gamma$ is bounded above by $2^{q(n)}$ for some polynomial $q(n)$. Finally, for every $\Gamma$ and $\varphi$, checking condition (\ref{Eq:Algorithm-Truth-Eq}) is an easy task. The analysis done so far shows that the time needed in the whole computation can be bounded above by $2^{2^{t(n)}}$ for some polynomial $t(n)$.
\end{proof}


A logic is said to be \textit{strongly finite} when it is determined by a finite set of finite matrices of finite type.\footnote{In some works the assumption of the finiteness of the type is dropped.} The next result provides upper bounds to the problem of classifying strongly finite logics in the Leibniz hierarchy of Figure \ref{Hierarchy}.

\begin{Corollary}
The problem of determining whether the logic of a finite set of finite matrices of finite type is algebraizable (resp.\ equivalential, protoalgebraic, weakly algebraizable) is in \textbf{EXPTIME}. The analogous problem for truth-equationality is in \textbf{2-EXPTIME}.
\end{Corollary}

\begin{proof}
The proof is a straightforward adaptation of those of Lemmas \ref{Lem:EXPTIME1} and \ref{Lem:EXPTIME2}. The only significant difference is that, in this case, the input is a set of matrices $\class{M}$ which are not necessarily reduced. For this reason, we need to replace the set $\class{M}$ with the set of reductions of its members. This amounts to compute some Leibniz congruences, and can be done efficiently (as explain in the proof of Lemma \ref{Lem:EXPTIME1}).
\end{proof}

\section{Two fundamental constructions}\label{Sec:Constructions}

In this section $\A$ is a fixed non-trivial algebra with basic operations $\mathcal{F}$, and $h$ is a unary function on $A$. For sake of simplicity, we assume that $\mathcal{F}$ contains no constant symbols. Our goal is to define two algebras $\A^{\natural}$ and $\A^{\flat}$ related to $\A$. The construction of these algebras presents some relation with the ones introduced in \cite{Ho91c} and \cite{FV09} to prove hardness results related to the study of type sets in tame congruence theory \cite{HoMcKe88} and to the study of Maltsev conditions \cite{Gratzer70Mal,Tay73}.

We begin by the definition of $\A^{\natural}$. The universe of $\A^{\natural}$ is given by eight disjoint copies $A_{1}, \dots, A_{8}$ of $A$. Given an element $a \in A$, we will denote by $a^{i}$ its copy in $A_{i}$. In this sense an arbitrary finite set of elements in $A^{\natural}$ can be denoted as $\{ a_{1}^{m_{1}}, \dots, a_{n}^{m_{n}} \}$ for some $a_{1}, \dots, a_{n} \in A$ and $m_{1}, \dots, m_{n} \leq 8$. The basic operations of $\A^{\natural}$ are the ones in $\mathcal{F}$ plus a new ternary operation $\heartsuit$ and a new unary operation $\Box$. Their interpretation is defined as follows. Given an $n$-ary operation $f \in \mathcal{F}$ and $a_{1}^{m_{1}} \dots, a_{n}^{m_{n}} \in A^{\natural}$, we set
\[
f(a_{1}^{m_{1}} \dots, a_{n}^{m_{n}}) \coloneqq f^{\A}(a_{1}, \dots, a_{n})^{5}.
\]
Observe that all the operations $f^{\A^{\natural}}$ with $f \in \mathcal{F}$ give values in $A_{5}$. Given $a^{m}, b^{n}, c^{k} \in A^{\natural}$, we set
\[
\heartsuit(a^{m}, b^{n}, c^{k})\coloneqq \left\{ \begin{array}{ll}
 a^{1} & \text{if $a^{m}=c^{k}$ and $h(a)^{5}= b^{n}$ and $m \in \{1, 3, 4\}$}\\
   a^{2} & \text{if $a^{m}=c^{k}$ and $h(a)^{5}= b^{n}$ and $m \in  \{ 2, 5, 6, 7, 8 \}$}\\
 a^{4} & \text{if $m, k \in \{1, 3, 4\}$ and (either $a^{m} \ne c^{k}$ or $h(a)^{5} \ne b^{n}$)}\\
 a^{7} & \text{if $\{ m, k \} \cap \{ 2, 5, 6, 7, 8 \} \ne \emptyset$ and}\\
 &\text{(either $a^{m} \ne c^{k}$ or $h(a)^{5} \ne b^{n}$).}\\
  \end{array} \right.  
\]
Given $a^{m} \in A^{\natural}$, we set
\[
\Box(a^{m})\coloneqq \left\{ \begin{array}{ll}
 a^{m} & \text{if $m=1$ or $m=2$}\\
 a^{m-1} & \text{if $m$ is even and $m \geq 3$}\\
  a^{m+1} & \text{if $m$ is odd and $m \geq 3$.}\\
  \end{array} \right.  
\]

Then we turn to define the algebra $\A^{\flat}$. The universe of $\A^{\flat}$ is given by two disjoint copies $A_{1}, A_{2}$ of $A$ plus two fresh elements $\{ 0, 1 \}$. The basic operation of $\A^{\flat}$ are the ones in $\mathcal{F}$ plus a new binary operation $+$, a new unary operation $\Box^{a}$ for every $a \in A$, and a new constant symbol $\boldsymbol{1}$. Their interpretation is defined as follows. The constant $\boldsymbol{1}$ denotes the element $1$. Given an $n$-ary operation $f \in \mathcal{F}$ and $b_{1}, \dots, b_{n} \in A^{\flat}$, we set
\[
f(b_{1}, \dots, b_{n})  \coloneqq \left\{ \begin{array}{ll}
 f^{\A}(a_{1}, \dots, a_{n})^{2} & \text{if $b_{i} \in \{ a_{i}^{1}, a_{i}^{2} \}$ for all $i \leq n$}\\
 1 & \text{if $b_{1} = \cdots = b_{n} = 1$}\\
 0 & \text{otherwise.}\\
  \end{array} \right.  
\]
Moreover, for every $a, b \in A^{\flat}$ we set
\[
a + b \coloneqq \left\{ \begin{array}{ll}
 1 & \text{either $a=b=1$ or}\\
 &\text{($a = c^{1}$ for some $c \in A$ and $b= h(c)^{2}$)}\\
 0 & \text{otherwise.}\\
  \end{array} \right.  
\]
Finally, for every $a \in A$ and $b \in A^{\flat}$ we set
\[
\Box^{a}(b) \coloneqq \left\{ \begin{array}{ll}
 b & \text{if $b \in (A_{1} \cup A_{2}) \smallsetminus \{ a^{1}, a^{2} \}$}\\
 a^{2} & \text{if $b =a^{1}$}\\
 a^{1} & \text{if $b =a^{2}$}\\
  0 & \text{if $b \in \{ 0, 1 \}$.}\\
  \end{array} \right.  
\]

Now, consider the matrices $\langle \A^{\natural}, F^{\natural}\rangle$ and $\langle \A^{\flat}, G^{\flat}\rangle$ where
\[
F^{\natural} \coloneqq A_{1} \cup A_{2} \text{ and }G^{\flat} \coloneqq A_{1} \cup \{ 1 \}.
\]

\begin{Lemma}\label{Lem:Reduced}
The matrices $\langle \A^{\natural}, F^{\natural}\rangle$ and $\langle \A^{\flat}, G^{\flat}\rangle$ are reduced.
\end{Lemma}

\begin{proof}
Observe that, to prove that an arbitrary matrix $\langle \B, H\rangle$ is reduced, it is enough to show that for every pair of different elements $a, b \in B$ there is a unary polynomial function $p(x)$ such that
\begin{equation}\label{Eq:Reduced}
\text{either }(p(a) \in H\text{ and }p(b) \notin H) \text{ or }(p(a) \notin  H\text{ and }p(b) \in  H).
\end{equation}

To prove that (\ref{Eq:Reduced}) holds for $\langle \A^{\natural}, F^{\natural}\rangle$ we reason as follows. Consider two distinct elements $a^{n}, b^{m} \in A^{\natural}$. Then consider the unary polynomial function
\[
p(x) \coloneqq \heartsuit(x, h(a)^{5}, a^{n}).
\]
It is clear that $p(a^{n}) \in F^{\natural}$, while $p(b^{m}) \notin F^{\natural}$. Hence $\langle \A^{\natural}, F^{\natural}\rangle$ is reduced.

Then we turn to prove that (\ref{Eq:Reduced}) holds for $\langle \A^{\flat}, G^{\flat}\rangle$ too. To this end, consider two different $a, b \in A^{\flat}$. If $a \in G^{\flat}$ and $b \notin G^{\flat}$, then we are done. Then suppose that either $a, b \in G^{\flat}$ or $a, b \notin G^{\flat}$. First consider the case where $a, b \notin G^{\flat}$. We can assume w.l.o.g.\ that $a = c^{2}$ for some $c \in A$. Then we have that $\Box^{c}(a) = c^{1} \in G^{\flat}$, while $\Box^{c}(b) = b \notin G^{\flat}$. Then suppose that $a, b \in G^{\flat}$. W.l.o.g.\ we have the following cases:
\benormal
\item $a= 1$ and $b \ne 1$.
\item $1 \notin \{ a, b \}$.
\enormal

1. Since $\A$ is non-trivial, there is $c \in A$ such that $b \ne c^{1}$. Then observe that $\Box^{c}(a) = 0 \notin G^{\flat}$, while $\Box^{c}(b) = b \in G^{\flat}$.

2. Then observe that $a=c^{1}$ and $b= d^{1}$ for some different $c, d \in A$. We have that $\Box^{c}(a) = c^{2} \notin G^{\flat}$, while $\Box^{c}(b) = b \in G^{\flat}$.
\end{proof}


\begin{Lemma}\label{Lem:Polynomial-Translation}
There is an algorithm that transform inputs $\langle \A, h \rangle$ of $\textsf{Gen-Clo}^{1}_{3}$ into matrices $\langle \A^{\natural}, F^{\natural}\rangle$ (resp.\ $\langle \A^{\flat}, G^{\flat}\rangle$) and runs in polynomial time in the length of the input.
\end{Lemma}

\begin{proof}
We detail the case of $\langle \A^{\natural}, F^{\natural}\rangle$ only. Observe that a typical input of $\textsf{Gen-Clo}^{1}_{3}$ is a pair $\langle \A, h \rangle$, where $\A$ is an algebra without constant symbols. Therefore it makes sense to transform it into a matrix of the form $\langle \A^{\natural}, F^{\natural}\rangle$.  Now, it is clear that this translation can be given an algorithmic form. Therefore, it only remains to prove that it runs in polynomial time. The only non-trivial part is that the length of the output $\langle \A^{\natural}, F^{\natural}\rangle$ should be bounded by a polynomial in the length $n$ of the input $\langle \A, h \rangle$.\  Since $h$ is unary, we know that $n$ dominates the cardinality of $\A$. Moreover,  $\A$ has at most $n$ basic operations, all of arity at most $3$. Thus $\A^{\natural}$ is an algebra of cardinality at most $8n$, with at most $n+2$ basic operations, all of arity at most $3$. Hence the length of $\A^{\natural}$ is bounded above by $8n + (n+2) \cdot (8n)^{3}$. Thus the length of $\langle \A^{\natural}, F^{\natural}\rangle$ is bounded by $8n + (n+2) \cdot (8n)^{3} + 8n$. 
\end{proof}

\section{Hardness results}

In this section we prove that the problem of determining whether the logic of a finite reduced matrix of finite type belongs to any level of the Leibniz hierarchy of Figure \ref{Hierarchy} is hard for \textbf{EXPTIME}. To this end, we will make use of the following purely algebraic decision problem related to clones \cite{Be11g,McMcTa87}:
\benbullet
\item $\textsf{Gen-clo}^{1}_{3}$: Given a finite non-trivial algebra $\A$ with finitely many basic operations $\mathcal{F}$, whose arities belong to the set $\{ 1, 2, 3 \}$, and a \textit{unary} operation $h$ on $A$ different from the identity, is $h$ in the clone of $\A$?
\ebullet
It is clear that $\textsf{Gen-clo}^{1}_{3}$ is decidable by means of a procedure which runs in exponential time. Remarkably, this procedure is the best possible \cite[Theorem 3.7]{BJS99}:\footnote{Observe that \cite[Theorem 3.7]{BJS99} states that a slightly different problem is complete for \textbf{EXPTIME}. But the proof given there shows that the same result holds for $\textsf{Gen-clo}^{1}_{3}$.}
\begin{Theorem}\label{Thm:Bergman}
The problem $\textsf{Gen-clo}^{1}_{3}$ is complete for \textbf{EXPTIME}.
\end{Theorem}


The proof of the next lemma relies on a series of technical results and one definition contained in the Appendix.

\begin{Lemma}\label{Lem:Polynomial}
Let $\A$ be a non-trivial algebra without constant symbols, $h$ be a unary operation on $A$ other than the identity, and $\vdash$ be the logic determined by the matrix $\langle \A^{\natural}, F^{\natural}\rangle$. The following conditions are equivalent:
\benroman
\item $\vdash$ is algebraizable.
\item $\vdash$ is weakly algebraizable.
\item $\vdash$ is equivalential.
\item $\vdash$ is protoalgebraic.
\item $h$ is in the clone of $\A$.
\eroman
\end{Lemma}

\begin{proof}
The implications (i)$\Rightarrow$(ii), (i)$\Rightarrow$(iii), (ii)$\Rightarrow$(iv) and (iii)$\Rightarrow$(iv) are straightforward. (v)$\Rightarrow$(i): Let $\varphi(x)$ be a term of $\A$ representing $h$ on $A$. Since $h$ is not the identity, we have that $\varphi(x) \ne x$. We consider the same term in the algebra $\A^{\natural}$. Then consider the following sets:
\[
\btau(x) \coloneqq \{ x \thickapprox \Box x \} \text{ and }\Delta(x, y) \coloneqq \{ \heartsuit(x, \varphi(x), y) \}.
\]
Bearing in mind that $\varphi(x) \ne x$, it is an easy exercise to check that for every $a, b \in A^{\natural}$ we have that
\begin{align*}
a \in F^{\natural} &\Longleftrightarrow \A^{\natural} \vDash \btau(a)\\
a = b &\Longleftrightarrow \Delta(a, b) \subseteq F^{\natural}.
\end{align*}
By condition 3 of Theorem \ref{Thm:Syntactic} we conclude that $\vdash$ is algebraizable.

(iv)$\Rightarrow$(v): Suppose that $\vdash$ is protoalgebraic. This means that there is a set $\Delta(x, y)$ satisfying the conditions in point 1 of Theorem \ref{Thm:Syntactic}. Since $\vdash$ is neither inconsistent nor almost-inconsistent, we know that $\Delta \ne \emptyset$ and that no variable belongs to $\Delta$. Pick any element $a \in A$. We know that $\{ a^{1} \} \cup \Delta(a^{1}, a^{1}) \subseteq F$, since $\emptyset \vdash \Delta(x, x)$. Together with the fact that $x, \Delta(x, y) \vdash y$ and that $a^{3} \notin F^{\natural}$, this implies that there is $\varphi(x, y) \in \Delta$ such that $\varphi(a^{1}, a^{1}) \in F^{\natural}$ and $\varphi(a^{1}, a^{3}) \notin F$.

Observe that $\emptyset \vdash \varphi(x, x)$. In particular, we have that $\varphi(a^{3}, a^{3}) \in F^{\natural}$. Together with the definition of $\A^{\natural}$, this implies that $\varphi(x, y) = \Box^{m}\heartsuit(\psi_{1}, \psi_{2}, \psi_{3})$ for some $m \in \omega$ and formulas $\psi_{1}, \psi_{2}$ and $\psi_{3}$. We can assume w.l.o.g.\ that $m=0$, since $\Box x \vdash x$ and $x \vdash \Box x$. The fact that $\varphi(x, x)$ is a tautology of $\vdash$ implies that
\begin{equation}\label{Eq:Condition}
\text{For every }b \in A^{\natural}\text{, if }c^{n} = \psi_{1}(b, b)\text{, then }\psi_{2}(b, b) = h(c)^{5}.
\end{equation}
Again looking at the definition of $\A^{\natural}$, this implies that $\psi_{2} = \Box^{2k}f(\epsilon_{1}, \dots, \epsilon_{n})$ for some $f \in \mathcal{F}$, formulas $\epsilon_{i}$ and $k \in \omega$. In particular, we can assume w.l.o.g.\ that $k=0$, since $\A^{\natural} \vDash  \Box^{2k}f(\epsilon_{1}, \dots, \epsilon_{n}) \thickapprox f(\epsilon_{1}, \dots, \epsilon_{n})$. To every formulas $\epsilon_{i}$ we associate a new formula $\delta_{i}$ written only with the symbols of $\mathcal{F}$. This is done as follows: first we delete from $\epsilon_{i}$ all occurrences of $\Box$, and then we replace every occurrence of $\heartsuit(\chi_{1}, \chi_{2}, \chi_{3})$ by $\chi_{1}$. This concludes the definition of $\delta_{i}$. Observe that for every $b, c \in A^{\natural}$ there is $a \in A$ such that
\[
\epsilon_{i}(b, c), \delta_{i}(b, c) \in \{ a^{1}, \dots, a^{8} \}.
\]
In particular, this implies that for every $b, c \in A^{\natural}$ we have
\[
f(\epsilon_{1}, \dots, \epsilon_{n})(b, c) = f(\delta_{1}, \dots, \delta_{n})(b, c).
\]
Thus we can assume w.l.o.g. that $\psi_{2} = f(\delta_{1}, \dots, \delta_{n})$.

Now, recall that $\varphi(a^{1}, a^{1}) \in F$, while $\varphi(a^{1}, a^{3}) \notin F$. In particular, we have that  $\varphi(a^{1}, a^{1}) \ne \varphi(a^{1}, a^{3})$. Then we can apply Lemma \ref{Lem:Central}, obtaining that $\varphi(a^{1}, a^{3}) \in A_{1} \cup A_{3} \cup A_{4}$. Together with Lemmas \ref{Lem:Tree} and \ref{Lem:Tree_Property} this implies that
\[
\varphi(b^{n}, b^{n}) \in \{ b^{1}, \dots, b^{8} \} \text{ for every }b^{n} \in A^{\natural}.
\]
Pick $b^{n} \in A^{\natural}$. In particular, we have that $\psi_{1}(b^{n}) = b^{m}$ for some $m \leq 8$. Together with condition (\ref{Eq:Condition}) this implies that 
\[
\psi_{2}(b^{n}, b^{n}) = f(\delta_{1}, \dots, \delta_{n})(b^{n}, b^{n}) = h(b)^{5}.
\]
Recall that $f(\delta_{1}, \dots, \delta_{n})$ is a formula written only with symbols in $\mathcal{F}$. Hence we obtain that
\[
f(\delta_{1}, \dots, \delta_{n})^{\A}(b, b) = h(b)
\]
 for every $b \in A$. This means that the term $f(\delta_{1}, \dots, \delta_{n})(x, x)$ represents $h$ in $\A$, i.e. $h$ belongs to the clone of $\A$.
\end{proof}

Let $\btau$ be a set of equations in variable $x$.\ A class of algebras $\class{K}$ is a $\btau$-\textit{algebraic semantics} for a logic $\vdash$ \cite{BP89,BlRe03,Ra06a} if for every set of formulas $\Gamma \cup \{ \varphi \}$, we have:
\[
\Gamma \vdash \varphi \Longleftrightarrow \btau(\Gamma) \vDash_{\class{K}} \btau(\varphi)
\]
where $\vDash_{\class{K}}$ is the equational consequence relative to $\class{K}$. In what follows we will make use of the following result \cite[Proposition 39]{Ra06a}:

\begin{Theorem}\label{Thm:JamesLemma}
Let $\vdash$ be a logic with a $\btau$-algebraic semantics. Suppose that $\btau$ contains an equation $\epsilon \thickapprox \delta$ such that $\delta(x) \vdash x$ and $\emptyset \vdash \epsilon(x)$. Then $\vdash$ is truth-equational by means of $\delta(x) \thickapprox \epsilon(x)$.
\end{Theorem}

\begin{Lemma}\label{Lem:Polynomial1}
Let $\A$ be a non-trivial algebra without constant symbols, $h$ be a unary operation of $A$, and $\vdash$ be the logic determined by the matrix $\langle \A^{\flat}, G^{\flat}\rangle$. $\vdash$ is truth-equational if and only if $h$ is in the clone of $\A$.
\end{Lemma}

\begin{proof}
We begin by the ``if'' part. Let $t(x)$ be the term that represents $h$ in $\A$. Then define $\delta \coloneqq x + t(x)$ and $\epsilon \coloneqq \boldsymbol{1}$. Moreover, set $\btau \coloneqq \{ \delta(x) \thickapprox \epsilon(x) \}$. Observe that the set of equations $\btau$ defines $G^{\flat}$ in $\A^{\flat}$. Thus the class $\class{K} = \{ \A^{\flat} \}$ is a $\btau$-algebraic semantics for $\vdash$. Moreover, observe that $\delta(x) \vdash x$ and $\emptyset \vdash \epsilon(x)$. Therefore we can apply \ref{Thm:JamesLemma} yielding that $\vdash$ is truth-equational.

Then we turn our attention to the ``only if'' part. Suppose that $\vdash$ is truth-equational. Observe that the matrix $\langle \A^{\flat}, G^{\flat}\rangle$ is reduced by Lemma \ref{Lem:Reduced}. Then we can use the technique described in the proof of Lemma \ref{Lem:EXPTIME2} to produce a set of equations $\btau$ that satisfies condition (\ref{Eq:TruthPredicates}), namely
\[
\btau \coloneqq \{ \epsilon(x) \thickapprox \delta(x) : \epsilon(a) = \delta(a) \text{ for every }a \in G^{\flat} \}.
\]

Since $\langle \A^{\flat}, G^{\flat}\rangle$ is reduced and $0 \in G^{\flat}$, there is an equation $\epsilon \thickapprox \delta \in \btau$ such that $\epsilon(0) \ne \delta(0)$. Looking at the definition of the basic operations of $\A^{\flat}$ we obtain that for every formula $\gamma(x)$, if $x$ really occurs in $\gamma$, then $\gamma(0) = 0$. Then either $x$ does not occur in $\epsilon$ or $x$ does not occur in $\delta$. We assume w.l.o.g.\ that $x$ does not occur in $\epsilon$. This means that $\epsilon$ is a closed term, i.e. a term of the form $t(\boldsymbol{1}, \dots, \boldsymbol{1})$ for some term $t(x_{1}, \dots, x_{n})$. In particular, we have that $\A^{\flat} \vDash t(\boldsymbol{1}, \dots, \boldsymbol{1}) \thickapprox \boldsymbol{1}$. Together with the definition of $\btau$, this implies that the equation $\delta(x) \thickapprox \boldsymbol{1}$ belongs to $\btau$.

Now, we claim that the symbol $\Box^{a}$ does not appear in $\delta(x)$, for every $a \in A$. To prove this, we will show by induction on the construction of $\delta$ that if the operation $\Box^{a}$ appears in $\delta$, then $\delta(1) = 0$. The case where $\delta$ is a variable is satisfied. Then suppose that $\delta = f(\varphi_{1}, \dots, \varphi_{n})$ for some $f \in \mathcal{F}$ and that $\Box^{a}$ appears in $\delta$. Then $\Box^{a}$ appears in some $\varphi_{i}$. By the inductive hypothesis we have that $\varphi_{i}(1) = 0$. Hence we conclude that
\[
\delta(1) = f(\varphi_{1}(a), \dots, \varphi_{i-1}(1), 0, \varphi_{i+1}(1), \dots, \varphi_{n}(1)) = 0.
\]
The case where $\delta = \varphi + \psi$ is handled analogously. Finally, consider the case where $\delta = \Box^{a} \varphi$. Observe that $\varphi(1) \in \{ 0, 1 \}$, since $\{ 0, 1 \}$ is the universe of a subalgebra of $\A^{\flat}$. Then $\Box^{a}(\varphi(1)) = 0$ as desired. This concludes the inductive proof. Now, recall that $1 \in G^{\flat}$ and, therefore, that $\delta(1) = 1$. Hence we conclude that $\Box^{a}$ does not appear in $\delta$, for every $a \in A$.

We are in the following situation: there is an equation $\delta \thickapprox \boldsymbol{1}$ in $\btau$ which satisfies the following requirements:
\benormal
\item $\Box^{a}$ does not appear in $\delta$.
\item $\delta$ is neither a closed term nor a variable.
\enormal 
The second condition above follows from the fact that $\delta(0) \ne 1$ and that $\delta(a) = 1$ for all $a \in A_{1}$. Conditions 1 and 2 imply that the principal symbol of $\delta$ must be either $f \in \mathcal{F}$ or $+$. If it is $f \in \mathcal{F}$, then $\delta = f(\varphi_{1}, \dots, \varphi_{n})$ for some formulas $\varphi_{1}, \dots, \varphi_{n}$. Looking at the definition of $f$ in $\A^{\flat}$, one sees that $\varphi_{i}(a) = 1$ for every $a \in G^{\flat}$. In particular, this means that all the equations $\varphi_{i} \thickapprox \boldsymbol{1}$ belong to $\btau$. Moreover, the fact that $\varphi_{i}(a) = 1$ for all $a \in A_{1}$ implies that $\varphi_{i}$ is not a variable. Finally, there must be at least one formula $\psi = \varphi_{i}$, which is not a closed term. Then this formula $\psi$ satisfies conditions 1 and 2 and, moreover, the equation $\psi \thickapprox \boldsymbol{1}$ belongs to $\btau$. Then consider the case where the principal symbol of $\delta$ is $+$, i.e. $\delta = \varphi + \psi$ for some terms $\varphi$ and $\psi$. Suppose that $\varphi$ is not a variable. The fact that the symbol $\Box^{a}$ does not appear in $\varphi$ for every $a \in A$ implies that $\varphi(a) \in A_{2} \cup \{ 0, 1 \}$ for every $a \in A^{\flat}$. Together with the fact that $\delta(a) = 1$ for every $a \in G^{\flat}$, this implies that $\varphi(a) = 1 = \psi(a)$ for every $a \in G^{\flat}$. This means that both the equations $\varphi \thickapprox \boldsymbol{1}$ and $\psi \thickapprox \boldsymbol{1}$ belong to $\btau$. Suppose that $\varphi$ does not satisfy condition 2. Since $\varphi$ cannot be a variable (since $\varphi(a) = 1$ for all $a \in A_{1}$), we conclude that $\varphi$ is the constant $\boldsymbol{1}$.  Since $\delta$ is not a closed term, this implies that $\psi$ is not a closed term. Moreover, $\psi$ cannot be a variable (since $\varphi(a) = 1$ for all $a \in A_{1}$). Hence we conclude that $\psi$ satisfies condition 2. Thus $\psi$ satisfies conditions 1 and 2 and, moreover, the equation $\psi \thickapprox \boldsymbol{1}$ belongs to $\btau$. Until now we have shown that either $\delta$ is of the form of $x + \psi$ or we can extract from it a proper subformula $\gamma$ which again satisfies conditions 1 and 2, and the equation $\gamma \thickapprox \boldsymbol{1}$ belongs to $\btau$. The same argument shows that either $\gamma$ is of the form $x + \psi$ or we can extract from it a proper subformula $\gamma'$, which satisfies conditions 1 and 2, and the equation $\gamma' \thickapprox \boldsymbol{1}$ belongs to $\btau$. Iterating this process we obtain a formula $x + \psi(x) \thickapprox \boldsymbol{1}$ which belongs $\btau$ and in which the symbols $\Box^{a}$ do not occur.

Consider any element $a^{1} \in A_{1}$. Since $a^{1} \in G^{\flat}$, we have that $a^{1} + \psi(a^{1}) = 1$. This means that $\psi(a^{1}) = h(a)^{2}$. In particular, observe that if the symbols $+$ or $\boldsymbol{1}$ appear in $\psi$, then we would have that $\psi(a^{1}) \in \{ 0, 1\}$, which is not the case. Thus the symbols $+$, $\boldsymbol{1}$ and $\Box^{a}$ do not occur in $\psi$. This means that $\psi$ is written with symbols among $\mathcal{F}$. Then for every $a \in A$ we have that
\[
\psi^{\A}(a)^{2} = \psi(a^{1}) = h(a)^{2}.
\]
Hence $\psi$ represents $h$ in $\A$, i.e. $h$ belongs to the clone of $\A$ as desired.
\end{proof}

We are now ready to state the main result of the paper:

\begin{Theorem}\label{Thm:Hard}
Let $\class{K}$ be a level in the Leibniz hierarchy of Figure \ref{Hierarchy}. The problem of determining whether the logic of a finite reduced matrix of finite type belongs to $\class{K}$ is hard for \textbf{EXPTIME}.
\end{Theorem}

\begin{proof}
Let $\class{K}$ be a level in the Leibniz hierarchy of Figure \ref{Hierarchy}. By Lemmas \ref{Lem:Polynomial} and \ref{Lem:Polynomial1} there is a many-one reduction of $\textsf{Gen-Clo}^{1}_{3}$ to the problem of determining whether the logic of a finite reduced matrix of finite type belongs to $\class{K}$. By Lemma \ref{Lem:Polynomial-Translation} this reduction runs in polynomial time. Together with Theorem \ref{Thm:Bergman}, this implies that the problem of determining whether the logic of a finite reduced matrix of finite type belongs to $\class{K}$ is hard for \textbf{EXPTIME}.
\end{proof}

\begin{Corollary}\label{Cor:EXPTIMEComplete}
Let $\class{K}$ be any level in the Leibniz hierarchy of Figure \ref{Hierarchy} different from the one of truth-equational logics. The problem of determining whether the logic of a finite reduced matrix of finite type belongs to $\class{K}$ is complete for \textbf{EXPTIME}.
\end{Corollary}

\begin{proof}
This is a combination of Lemma \ref{Lem:EXPTIME1} and Theorem \ref{Thm:Hard}.
\end{proof}

Remarkably, the methods developed in this work can be adapted to prove hardness results for other classes of logics studied in abstract algebraic logic. To exemplify this phenomenon, recall that a logic $\vdash$ is \textit{order algebraizable} \cite{Ra06,JRa13} if there are a set $\Delta(x, y)$ of formulas  and a set $\btau(x)$ of inequalities such that for every $\langle \A, F\rangle \in \Modstar(\vdash)$ the relation $\leq_{\A}^{F}$ defined for every $a, b \in A$ as
\[
a \leq_{\A}^{F} b \Longleftrightarrow \Delta(a, b) \subseteq F
\]
is a partial order, and for every $a \in A$
\[
a \in F \Longleftrightarrow \langle \A, \leq_{\A}^{F}\rangle \vDash \btau(a).
\]
\begin{Corollary}
The problem of determining whether the logic of a finite reduced matrix of finite type is order algebraizable is complete for \textbf{EXPTIME}.
\end{Corollary}

\begin{proof}[Sketch of the proof.]
The class of order algebraizable logics has a syntactic characterization \cite[Theorem 7.1]{JRa13}. Applying an argument similar to the one described in the proof of Lemma \ref{Lem:EXPTIME1} to this characterization, it is possible to show that the problem of determining whether the logic of a finite reduced matrix of finite type is order algebraizable belongs to \textbf{EXPTIME}. To prove that this problem is hard for \textbf{EXPTIME}, we reason as follows. First observe that every order algebraizable logic is is equivalential, and that every algebraizable logic is order algebraizable. Then we can add the fact that $\vdash$ is order algebraizable to the list of equivalent conditions in the statement of Lemma \ref{Lem:Polynomial}. This yields the desired hardness result.
\end{proof}

We conclude with a list of open problems on the decidability and computational complexity of the hierarchies of abstract algebraic logic.

\begin{problem}
Is the problem of determining whether the logic of a finite reduced matrix of finite type is truth-equational in \textbf{EXPTIME}? Observe Theorem \ref{Thm:Hard} shows that this problem is hard for \textbf{EXPTIME}.
\end{problem}

\begin{problem}
Is the problem of determining whether the logic of a finite reduced matrix of finite type has an \textit{algebraic semantics} decidable? And if so, then which is its computational complexity?
\end{problem}

\begin{problem}
The problem of determining whether the logic of a finite reduced matrix of finite type is \textit{assertional} \cite{AFRM15} is easily seen to be decidable. More interestingly, which is its computational complexity? Similar questions can be risen about other classes of logics in the Leibniz hierarchy.
\end{problem}

\begin{problem}
In this work we focussed on logics determined by arbitrary finite reduced matrices of finite type. Are the hardness results still true if we restrict our attention to finite reduced matrices $\langle \A, F\rangle$ of finite type such that $F$ is a \textit{singleton}?
\end{problem}

\begin{problem}
Are there reasonable assumptions (on matrices) under which the problem of classifying the logic of a finite reduced matrix of finite type in the Leibniz hierarchy becomes \textit{tractable}?
\end{problem}

\begin{problem}
Investigate the decidability and complexity of determining whether the logic of a finite reduced matrix of finite type belongs to a given level of the Frege hierarchy of abstract algebraic logic.
\end{problem}

\paragraph{\bfseries Acknowledgements.}
Thanks are due to Carles Noguera for rising the question about the computational complexity of the problem of classifying logics in the Leibniz hierarchy, and for providing very useful comments on early versions of the paper, which contributed to improve the presentation. We also also very grateful to the anonymous referee for a variety of helpful observations, including a correction to point 3 of Theorem 2.1. The author was supported by  the grant GBP$202/12/$G$061$ of the Czech Science Foundation, and by project CZ.$02$.$2$.$69$/$0$.$0$/$0$.$0$/$17\_050$/$0008361$, OPVVV M\v{S}MT, MSCA-IF Lidske\'e zdroje v teoreticke\'e informatice.
\section{Appendix}

Through this section we assume that $\A$ is a fixed non-trivial algebra without constant symbols, and $\A^{\natural}$ is the algebra associated with it as in Section \ref{Sec:Constructions}. The \textit{subformula tree} of a formula $\varphi$ of $\A^{\natural}$ is defined by recursion on the construction of $\varphi$ as follows. The subformula tree of a variable $x$ is the one-element tree, whose unique node is labelled by $x$. Then the subformula tree of a complex formula $\varphi$ is obtained as follows. Suppose that $\varphi = g(\psi_{1}, \dots, \psi_{n})$ for some basic $n$-ary symbol $g$ and formulas $\psi_{1}, \dots, \psi_{n}$. First we pick the disjoint union of the subformula trees of $\psi_{1}, \dots, \psi_{n}$, and we relabel the root of the subformula tree of $\psi_{i}$ by $\langle \psi_{i}, i \rangle$. Second we add to these trees a common root labelled by $\varphi$.

\begin{law}
A formula $\varphi$ of $\A^{\natural}$ has the \textit{tree property} when every node of the subformula tree of $\varphi$, which is labelled by a formula whose principal symbol belongs to $\mathcal{F}$ (plus possibly a natural number), is preceded or equal to a point which has a label of the form $\langle \beta, 2\rangle$ whose immediate predecessor is labelled by $\heartsuit(\alpha, \beta, \gamma)$ or $\langle \heartsuit(\alpha, \beta, \gamma), n\rangle$ for some $n \in \omega$.
\end{law}

\begin{exa}
Since the definition of tree property is quite technical, we proceed to explain it through some basic examples. For the sake of simplicity, suppose that $\mathcal{F}$ contains only a unary symbol $f(x)$. The formulas
\[
\heartsuit(x, \Box f(x), y) \text{ and } \heartsuit( \Box z, \heartsuit(f(z),  f(f(x)), y), z)
\]
have the tree property. On the other hand, the formula $\heartsuit(f(x), f(x), z)$ does not have the tree property. To explain why, observe that its subformula tree looks as follows:
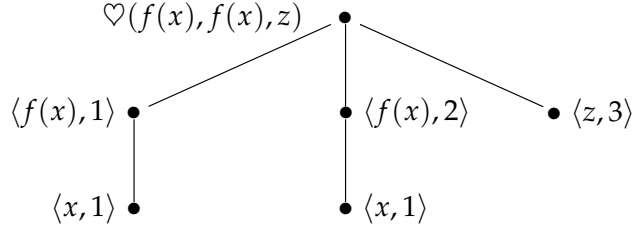
\[
\xymatrix@R=36pt @C=80pt @!0{
&*-{\heartsuit(f(x), f(x), z) \ \ \ \  \bullet\  \ \ \ \ \ \ \ \ \ \ \ \ \ \ \ \ \ \ \ \ \  \ \ \ \ \ \   }\ar@{-}[d]\ar@{-}[dr] \ar@{-}[dl]&\\
*-{\langle f(x), 1\rangle  \ \bullet  \ \ \ \ \ \ \ \ \ \ \ \ \ \   }\ar@{-}[d]&*-{ \  \ \ \ \ \ \ \ \ \ \ \ \ \ \bullet\  \langle f(x), 2\rangle  }\ar@{-}[d]&*-{ \ \ \ \ \ \ \ \ \bullet \ \langle z, 3\rangle   }\\
*-{\langle x, 1\rangle  \ \bullet   \ \ \ \ \ \ \ \ \ }&*-{ \ \ \ \ \ \ \ \ \ \bullet \ \langle x, 1\rangle   }&
}
\]
It is not hard to see that the occurrence of $\langle f(x), 1 \rangle$ in the lefter branch makes the tree property fail.
\qed 
\end{exa}

\begin{Lemma}\label{Lem:Tree_Property}
If $\varphi(x_{1}, \dots, x_{n})$ has the tree property, then for every $b^{n} \in A^{\natural}$ there is $m \leq 8$ such that
\[
\varphi(b^{n}, \dots, b^{n}) = b^{m}.
\]
Moreover, if $i \in \{1, 3, 4\}$, then $\varphi(b^{i}, \dots, b^{i}) \in \{ b^{1}, b^{3}, b^{4} \}$.
\end{Lemma}

\begin{proof}
We reason by induction on the construction of $\varphi$. If $\varphi$ is a variable, then the result holds vacuously. Then consider the case where $\varphi$ is a complex formula. If $\varphi = f(\psi_{1}, \dots, \psi_{n})$ for some $f \in \mathcal{F}$, then clearly $\varphi$ does not have the tree property. Then consider the case where $\varphi = \Box \psi$. Then clearly also $\psi$ has the tree property. Consider $b^{n} \in A^{\natural}$. By inductive hypothesis there is $m \leq 8$ such that $\psi(b^{n}, \dots, b^{n}) = b^{m}$. Together with the definition of $\Box$, this implies that
\[
\varphi(b^{n}, \dots, b^{n}) = \Box\psi(b^{n}, \dots, b^{n}) = \Box b^{m} \in \{ b^{1}, \dots, b^{8} \}.
\]
Finally, consider the case where $\varphi = \heartsuit(\psi_{1}, \psi_{2}, \psi_{3})$. Clearly $\psi_{1}$ and $\psi_{3}$ must have the tree property. We focus only  on $\psi_{1}$. Consider $b^{n} \in A^{\natural}$. By inductive hypothesis there is $m \leq 8$ such that $\psi_{1}(b^{n}, \dots, b^{n}) = b^{m}$. Together with the definition of $\heartsuit$, this implies that
\[
\heartsuit(\psi_{1}, \psi_{2}, \psi_{3})(b^{n}, \dots, b^{n}) \in \{ b^{1}, \dots, b^{8}\}.
\]
This concludes the inductive argument.

The last part of the statement (the one starting with \textit{Moreover...}) is also proved by a simple induction.
\end{proof}

\begin{Lemma}\label{Lem:Tree}
Let $\varphi(x, y)$ be a formula and $a \in A$. $\varphi(a^{1}, a^{3}) \in A_{1} \cup A_{3} \cup A_{4}$ if and only if $\varphi$ has the tree property.
\end{Lemma}

\begin{proof}
The direction from left to right is proved by induction on the construction of $\varphi$. The case where $\varphi$ is a variable is clear. Then suppose that $\varphi$ is a complex formula. If $\varphi = f(\psi_{1}, \dots, \psi_{n})$ for some $f \in \mathcal{F}$, then $\varphi(a^{1}, a^{3})$ belongs to $A_{5}$, which is disjoint from $A_{1} \cup A_{3} \cup A_{4}$. Then suppose that $\varphi = \Box \psi$. From $\Box \psi(a^{1}, a^{3}) \in A_{1} \cup A_{3} \cup A_{4}$ and the definition of $\Box$, we can infer that $\psi(a^{1}, a^{3}) \in A_{1} \cup A_{3} \cup A_{4}$. By the inductive hypothesis we obtain that $\psi$ has the tree property. In particular, this is implies that $\varphi$ has the tree property too. Then consider the case where $\varphi = \heartsuit(\psi_{1}, \psi_{2}, \psi_{3})$. From the fact that $\varphi(a^{1}, a^{3}) \in A_{1} \cup A_{3} \cup A_{4}$, we can infer that $\psi_{1}(a^{1}, a^{3}), \psi_{3}(a^{1}, a^{3}) \in A_{1} \cup A_{3} \cup A_{4}$. By the inductive hypothesis, this means that both $\psi_{1}$ and $\psi_{3}$ have the tree property. In particular, this implies that $\heartsuit(\psi_{1}, \psi_{2}, \psi_{3})$ has the tree property too.

The direction from right to left also follows from a simple induction.
\end{proof}

\begin{Lemma}\label{Lem:Easy_fact}
For every formula $\psi (x, \vec{y})$, $\vec{c} \in A^{\natural}$, and $a \in A$, there is $b \in A$ and $n, m \leq 8$ such that $\psi(a^{1}, \vec{c}) = b^{n}$ and $\psi(a^{3}, \vec{c}) = b^{m}$.
\end{Lemma}

\begin{proof}
This follows from a simple induction on the construction of $\psi$, using the definition of the basic operations of $\A^{\natural}$.
\end{proof}

As a consequence we obtain the following:

\begin{Corollary}\label{Cor:Equal}
Let $f \in \mathcal{F}$ be $n$-ary and $\varphi_{1}(x, \vec{y}), \dots, \varphi_{n}(x, \vec{y})$ be terms of $\A^{\natural}$. For every $a \in A$ and $\vec{c} \in A^{\natural}$ we have
\[
f(\varphi_{1}, \dots, \varphi_{n})(a^{1}, \vec{c}) = f(\varphi_{1}, \dots, \varphi_{n})(a^{3}, \vec{c}).
\]
\end{Corollary}
\begin{proof}
By Lemma \ref{Lem:Easy_fact} there are $b_{1},\dots, b_{n} \in A$ and $k_{1}, \dots, k_{n}, m_{1}, \dots, m_{n} \leq 8$ such that
\[
\varphi_{i}(a^{1}, \vec{c}) = b_{i}^{m_{i}} \text{ and }\varphi_{i}(a^{3}, \vec{c}) = b_{i}^{k_{i}}.
\]
This implies that
\begin{align*}
f(\varphi_{1}, \dots, \varphi_{n})(a^{1}, \vec{c}) &= f(b_{1}^{m_{1}}, \dots, b_{n}^{m_{n}})\\
&= f^{\A}(b_{1}, \dots, b_{n})^{5}\\
&=f(b_{1}^{k_{1}}, \dots, b_{n}^{k_{n}})\\
&= f(\varphi_{1}, \dots, \varphi_{n})(a^{3}, \vec{c})
\end{align*}
concluding the proof.
\end{proof}

\begin{Lemma}\label{Lem:Central}
Let $\varphi(x, y)$ be a formula and $a \in A$. If $\varphi(a^{1}, a^{3}) \ne \varphi(a^{1}, a^{1})$, then $\varphi(a^{1}, a^{3}) \in A_{1} \cup A_{3} \cup A_{4}$.
\end{Lemma}

\begin{proof}
We reason by induction on the construction of $\varphi$. The case where it is a variable is  satisfied. Then suppose that $\varphi$ is a complex formula. If $\varphi = f(\psi_{1}, \dots, \psi_{n})$ for some $f \in \mathcal{F}$, then by Corollary \ref{Cor:Equal} we obtain that $\varphi(a^{1}, a^{3}) = \varphi(a^{1}, a^{1})$, which is false. Then suppose that $\varphi = \Box \psi$. From $\Box \psi(a^{1}, a^{3})\ne \psi(a^{1}, a^{1})$ it follows that  $\psi(a^{1}, a^{3}) \ne \psi(a^{1}, a^{1})$. By inductive hypothesis we obtain that $\psi(a^{1}, a^{3}) \in A_{1} \cup A_{3} \cup A_{4}$. Together with the definition of $\Box$, this implies that $\varphi(a^{1}, a^{3}) \in A_{1} \cup A_{3} \cup A_{4}$ as well. It only remains to consider the case where $\varphi = \heartsuit(\psi_{1}, \psi_{2}, \psi_{3})$. From the fact that $\varphi(a^{1}, a^{3}) \ne \varphi(a^{1}, a^{1})$ it follows that one of the following conditions holds:
\benormal
\item $\psi_{1}(a^{1}, a^{3}) \ne \psi_{1}(a^{1}, a^{1})$.
\item $\psi_{2}(a^{1}, a^{3}) \ne \psi_{2}(a^{1}, a^{1})$ and 1 does not hold.
\item $\psi_{3}(a^{1}, a^{3}) \ne \psi_{3}(a^{1}, a^{1})$ and neither 1 not 2 hold.
\enormal
We consider these cases separately:

1. By the inductive hypothesis we have that $\psi_{1}(a^{1}, a^{3}) \in A_{1} \cup A_{3} \cup A_{4}$. If $\psi_{3}(a^{1}, a^{3}) \in A_{1} \cup A_{3} \cup A_{4}$, then by definition of $\heartsuit$ we conclude that $\varphi(a^{1}, a^{3}) \in A_{1} \cup A_{3} \cup A_{4}$ and we are done. Then suppose towards a contradiction that $\psi_{3}(a^{1}, a^{3}) \notin  A_{1} \cup A_{3} \cup A_{4}$. By Lemma \ref{Lem:Tree} we know that $\psi_{3}$ does not have the tree property. Together with the definition o $\A^{\natural}$, this implies that $\psi_{3}(a^{1}, a^{1}) \notin A_{1} \cup A_{3} \cup A_{4}$ too.  Moreover, by Lemma \ref{Lem:Tree} we know that $\psi_{1}$ has the tree property. Together with Lemma \ref{Lem:Tree_Property} this implies that
\[
\psi_{1}(a^{1}, a^{1}), \psi_{1}(a^{1}, a^{3}) \in \{ a^{1}, a^{3}, a^{4} \}.
\]
Thus we obtain that
\[
\heartsuit(\psi_{1}, \psi_{2}, \psi_{3})(a^{1}, a^{3}) = a^{7} = \heartsuit(\psi_{1}, \psi_{2}, \psi_{3})(a^{1}, a^{1})
\]
contradicting the fact that $\varphi(a^{1}, a^{3}) \ne \varphi(a^{1}, a^{1})$.

2. We will show that this case leads to a contradiction, that is to say that it never happens. If $\psi_{2}(a^{1}, a^{3}) \ne \psi_{2}(a^{1}, a^{1})$, then by inductive hypothesis this means that $\psi_{2}(a^{1}, a^{3}) \in A_{1} \cup A_{3} \cup A_{4}$. By Lemma \ref{Lem:Tree} we know that $\psi_{2}$ has the tree property and by Lemma \ref{Lem:Tree_Property} that $\psi_{2}(a^{1}, a^{1}) \in A_{1} \cup A_{3} \cup A_{4}$ as well. Consider $b \in A$ and $n \leq 8$ such that $\psi_{1}(a^{1}, a^{1}) = \psi_{1}(a^{1}, c^{3}) = b^{n}$. By definition of $\heartsuit$ we have that
\[
\varphi(a^{1}, a^{1}) =  \left\{ \begin{array}{ll}
b^{4} & \text{if $\psi_{1}(a^{1}, a^{1}) , \psi_{3}(a^{1}, a^{1}) \in A_{1} \cup A_{3} \cup A_{4}$}\\
 b^{7} & \text{otherwise.}\\
  \end{array} \right.  
  \]
  \[
  \varphi(a^{1}, a^{3}) =  \left\{ \begin{array}{ll}
b^{4} & \text{if $\psi_{1}(a^{1}, a^{3}) , \psi_{3}(a^{1}, a^{3}) \in A_{1} \cup A_{3} \cup A_{4}$}\\
 b^{7} & \text{otherwise.}\\
  \end{array} \right.  
\]
Since $\psi_{1}(a^{1}, a^{1}) = \psi_{1}(a^{1}, a^{3})$ by assumption and $\varphi(a^{1}, a^{1}) \ne  \varphi(a^{1}, a^{3})$, we have that $\psi_{3}(a^{1}, a^{1}) \ne \psi_{3}(a^{1}, a^{3})$. By inductive hypothesis this means that $\psi_{3}(a^{1}, a^{3}) \in A_{1} \cup A_{3} \cup A_{4}$. By Lemmas \ref{Lem:Tree} and \ref{Lem:Tree} we obtain that $\psi_{3}(a^{1}, a^{1}) \in A_{1} \cup A_{3} \cup A_{4}$ too. But this implies that 
\[
\varphi(a^{1}, a^{1}) = a^{4} = \varphi(a^{1}, a^{3})
\]
contradicting the fact that $\varphi(a^{1}, a^{1}) \ne \varphi(a^{1}, a^{3})$.

3. By inductive hypothesis we have that $\psi_{3}(a^{1}, a^{3}) \in A_{1} \cup A_{3} \cup A_{4}$. Suppose towards a contradiction that $\varphi(a^{1}, a^{3}) \notin A_{1} \cup A_{3} \cup A_{4}$. Then we have that $\psi_{1}(a^{1}, a^{3}) \notin A_{1} \cup A_{3} \cup A_{4}$. Consider $b \in A$ and $n \leq 8$ such that $b^{n} = \psi_{1}(a^{1}, a^{3})$. By the definition of $\heartsuit$ we have that $\varphi(a^{1}, a^{3}) = b^{7}$. On the other hand, observe that applying Lemmas \ref{Lem:Tree} and \ref{Lem:Tree_Property} to the fact that $\psi_{3}(a^{1}, a^{3}) \in A_{1} \cup A_{3} \cup A_{4}$, we obtain that $\psi_{3}(a^{1}, a^{1}) \in A_{1} \cup A_{3} \cup A_{4}$ as well. Together with the fact that $\psi_{1}(a^{1}, a^{1}) = \psi_{1}(a^{1}, a^{3}) \notin A_{1} \cup A_{3} \cup A_{4}$, this implies that $\varphi(a^{1}, a^{1}) = b^{7}$. But this contradicts the fact that $\varphi(a^{1}, a^{1}) \ne \varphi(a^{1}, a^{3})$.
\end{proof}

\bibliographystyle{plain}

\end{document}